\documentclass[a4paper]{article}

\usepackage{enumitem}
\usepackage[utf8]{inputenc}
\usepackage[T1]{fontenc}
\usepackage{textcomp}
\usepackage{lmodern}
\usepackage{babel}
\usepackage{lipsum}
\usepackage{graphicx}
\usepackage{scalerel,stackengine,amsmath}
\usepackage{subcaption}
\usepackage{tikz}
\usetikzlibrary{positioning, calc}
\usepackage{algorithm}
\usepackage{algpseudocode}
\usepackage{accents}
\usepackage{xcolor,color}
\usepackage{parskip}
\usepackage{bm}
\usepackage{amssymb,amsmath,amsthm,amsfonts,mathtools}
\usepackage[hyperfootnotes=false]{hyperref}
\usepackage{thmtools}
\usepackage{stmaryrd}
\usepackage{multirow}
\usepackage{booktabs}
\usepackage{alphabeta} 
\usepackage{thm-restate}
\usepackage{nicefrac}

\newif\ifarxiv
\arxivfalse

\newif\ifopus
\opustrue

\ifarxiv
    \newif\ifarxivoropus
        \arxivoropustrue
\else
    \ifopus
        \newif\ifarxivoropus
            \arxivoropustrue
    \else
        \newif\ifarxivoropus
            \arxivoropusfalse
    \fi
\fi

\usepackage[scaled=0.95]{FiraMono}
\ifarxiv
\usepackage{arxiv}
\usepackage{charter}
\usepackage{fullpage}
\fi

\ifopus
\usepackage{zibtitlepage}
\usepackage{arxiv}
\usepackage{charter}
\usepackage{fullpage}
\fi

\hypersetup{
    colorlinks=true,
    linkcolor=blue,
    filecolor=magenta,
    urlcolor=cyan,
    citecolor=blue
    }

\ifarxivoropus
\author{
    \href{https://orcid.org/0000-0003-3673-966X}{\includegraphics[scale=0.06]{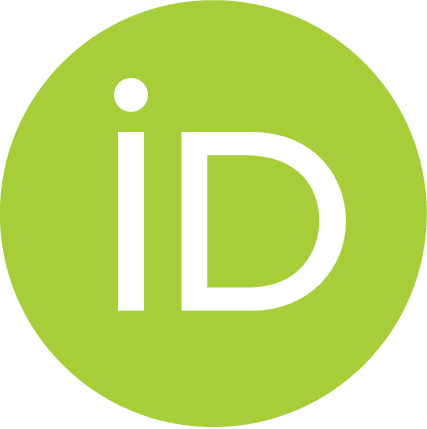}}\hspace{1mm}Gabriele Iommazzo\thanks{LIX CNRS, \'Ecole Polytechnique, Institut Polytechnique de Paris, Palaiseau, France. Gabriele Iommazzo is now at Zuse Institute Berlin, Berlin, Germany (\texttt{iommazzo@zib.de})}\\
    \And
    \href{https://orcid.org/0000-0002-4040-0960}{\includegraphics[scale=0.06]{orcid_id_icon.eps}}\hspace{1mm}Claudia D'Ambrosio\footnotemark[2]\\
    \texttt{dambrosio@lix.polytechnique.fr} \\
    \And
    \href{https://orcid.org/0000-0002-5704-3170}{\includegraphics[scale=0.06]{orcid_id_icon.eps}}\hspace{1mm}Antonio Frangioni\thanks{Dip.~di Informatica, Universit\`a di Pisa, Pisa, Italy} \\
    \texttt{frangio@di.unipi.it} \\
    \And
    \href{https://orcid.org/0000-0003-3139-6821}{\includegraphics[scale=0.06]{orcid_id_icon.eps}}\hspace{1mm}Leo Liberti\footnotemark[2]\\
    \texttt{liberti@lix.polytechnique.fr} \\
}
\else
\author{Gabriele Iommazzo\inst{1}\orcidID{0000-0003-3673-966X} \and
Claudia D'Ambrosio\inst{2}\orcidID{0000-0002-4040-0960}\and
Antonio Frangioni\inst{3}\orcidID{0000-0002-5704-3170}\and
Leo Liberti\inst{2}\orcidID{0000-0003-3139-6821}
}
\authorrunning{G. Iommazzo et al.}
\institute{LIX CNRS, \'Ecole Polytechnique, Institut Polytechnique de Paris, Palaiseau, France\and
Dip.~di Informatica, Università di Pisa, Pisa, Italy\\
\email{iommazzo@zib.de, \{dambrosio,giommazz,liberti\}@lix.polytechnique.fr}}
}
\fi

\ifopus
\ZTPAuthor{
	\ZTPHasOrcid{Gabriele Iommazzo}{0000-0003-3673-966X},
	\ZTPHasOrcid{Claudia D'Ambrosio}{0000-0002-4040-0960},
	\ZTPHasOrcid{Antonio Frangioni}{0000-0002-5704-3170},
	\ZTPHasOrcid{Leo Liberti}{0000-0003-3139-6821},}
\ZTPTitle{A learning-based mathematical programming formulation for the automatic configuration of optimization solvers}
\ZTPNumber{23-28}
\ZTPMonth{December}
\ZTPYear{2023}
\fi

\title{A learning-based mathematical programming formulation for the automatic configuration of optimization solvers
\thanks{This paper has received funding from the European Union's Horizon 2020 research and innovation programme under the Marie Sklodowska-Curie grant agreement n.~764759 ``MINOA''.}}

\begin{document}


\maketitle

\vspace{5mm}

\begin{center}
\begin{minipage}{0.85\textwidth}
\begin{center}
 \textbf{Abstract}
\end{center}
 {\small We propose a methodology, based on machine learning and optimization, for selecting a solver configuration for a given instance. First, we employ a set of solved instances and configurations in order to learn a performance function of the solver. Secondly, we formulate a mixed-integer nonlinear program where the objective/constraints explicitly encode the learnt information, and which we solve, upon the arrival of an unknown instance, to find the best solver configuration for that instance, based on the performance function. The main novelty of our approach lies in the fact that the configuration set search problem is formulated as a mathematical program, which allows us to a) enforce hard dependence and compatibility constraints on the configurations, and b) solve it efficiently with off-the-shelf optimization tools.
}
\end{minipage}
\end{center}

\section{Introduction}\label{s:intro}

We address the problem of finding instance-wise optimal configurations for general Mathematical Programming (MP) solvers. We are particularly motivated by state-of-the-art general-purpose solvers, which combine a large set of diverse algorithmic components (relaxations, heuristics, cutting planes, branching, \ldots) and therefore have a long list of user-configurable parameters; tweaking them can have a significant impact on the quality of the obtained solution and/or on the efficiency of the solution process (see,~e.g., \cite{HHL10}). Good solvers have effective default parameter configurations, carefully selected to provide good performances in most cases. Furthermore, tuning tools may be available (e.g., \cite[Ch.~10]{cplex127}) which run the solver, with different configurations, on one or more instances within a given time limit, and record the best parameter values encountered. Despite all this, the produced parameter configurations may still be highly suboptimal with specific instances. Hence, a manual search for the best parameter values may be required. This is a highly nontrivial and time-consuming task, due to the large amount of available parameters (see, e.g., \cite{cplex_param_ref}), which requires a profound knowledge of the application at hand and an extensive experience in solver usage. Therefore, it is of significant interest to develop general approaches, capable of performing it efficiently and effectively in an automatic way.

This setting is an instance of the Algorithm  Configuration Problem (ACP)  \cite{rice_acp}, which is defined as follows: given a target algorithm, its set of parameters, a set of instances of a problem class and a measure of the performance of the target algorithm on a pair (instance, algorithmic configuration), find the parameter configuration providing optimal algorithmic performance according to the given measure, on a specific instance or instance set. Several domains can benefit from automating this task. Some possible applications, beyond MP solvers (see, e.g., \cite{ac_best_practices} and references therein), are: solver configuration for the propositional satisfiability problem, hyperparameter tuning of ML models or pipelines, algorithm selection, administering ad-hoc medical treatment, etc.
Our approach for addressing the ACP on MP solvers is based on a two-fold process:
\begin{enumerate}[label=(\roman*)]
	\item in the \emph{Performance Map Learning Phase} (PMLP), supervised Machine Learning (ML) techniques \cite{MRT18} are used to learn a \emph{performance function}, which maps some features of the instance being solved, and the parameter configuration, into some measure of solver efficiency and effectiveness;
	\item the formal model underlying the ML methodology used in the PMLP is translated into MP terms; the resulting formulation, together with constraints encoding the compatibility of the configuration parameter values, yields the \emph{Configuration Set Search Problem} (CSSP), a Mixed-Integer Nonlinear Program (MINLP) which, for a given instance, finds the configuration providing optimal performance with respect to the performance function.
\end{enumerate}
The main novelty of our approach lies in the fact that we explicitly model and optimize the CSSP using the mathematical description of the PMLP technique. This is in contrast to most of the existing ACP approaches, which instead employ heuristics such as local searches \cite{calibra,paramils}, genetic algorithms \cite{gga}, evolutionary strategies \cite{lao} and other methods \cite{irace}. Basically, most approaches consider the performance function as a black box, even when it is estimated by means of some ML technique and, therefore, they cannot reasonably hope to find a global minimum when the number of parameters grows. Rather, one of the strengths of our methodology is that it exploits the mathematical structure of the CSSP, solving it with sophisticated, off-the-shelf MP solvers. Moreover, formulating the CSSP by MP is advantageous as it allows the seamless integration of the compatibility constraints on the configuration parameters, which is something that other ACP methods may struggle with.
The idea of using a ML predictor to define the unknown components (constraints, objective) of a MP has been already explored in \emph{data-driven optimization}.
In general, it is possible to represent the ML model of a mapping/relation as a MP (or, equivalently, a Constraint Programming model) and optimize upon this \cite{edml}; however, while this is in principle possible, the set of successful applications in practice is limited. Indeed, using this approach in the ACP context is, to the best of our knowledge, new; it also comes with some specific twists.
We tested this idea with the following components: we configured nine parameters of the IBM ILOG CPLEX solver \cite{cplex127}, which we employed to solve instances of the Hydro Unit Commitment (HUC) problem \cite{huc}, we chose Support Vector Regression (SVR) \cite{svr} as the PMLP learning methodology, and we used the off-the-shelf MINLP solver Bonmin \cite{bonmin_algo} to solve the CSSP.

The paper is structured as follows: in Sec.~\ref{s:acp} we will review existing work on algorithm configuration; in Sec.~\ref{s:pmlp_cssp} we will detail our approach and provide the explicit formulation of the CSSP with SVR; in Sec.~\ref{s:experiments} we will discuss some computational results.

\section{The algorithm configuration problem}\label{s:acp}

Most algorithms have a very high number of configurable parameters of various types (boolean, categorical, integer, continuous), which usually makes the ACP very hard to solve in practice. Notably, this issue significantly affects MP solvers: they are highly complex pieces of software, embedding several computational components that tackle the different phases of the solution process; the many available algorithmic choices are exposed to the user as a long list of tunable parameters (for example,  more than 150 in CPLEX \cite{cplex_param_ref}).

\smallskip
Approaches to the ACP  can be compared based on how they fit into the following two categories: Per-Set (PS) or Per-Instance (PI); offline or online.
In PS approaches, the optimal configuration is defined as the one with the best overall performance over a set of instances belonging to the same problem class. Therefore, PS approaches first find the optimal configuration for a problem class and then use it for any instance pertaining to that class. The exploration of the configuration set is generally conducted by means of heuristics, such as various local search procedures \cite{calibra,paramils}, racing methods \cite{irace}, genetic algorithms \cite{gga} or other evolutionary algorithms \cite{revac}. In this context, an exception is, e.g., the approached described in \cite{smac}, which predicts the performance of the target algorithm by random forest regression and then uses it to guide the sampling in an iterative local search.
PS approaches, however, struggle when the target algorithm performance varies considerably among instances belonging to the same problem class. In these cases, PI methodologies, which assume that the optimal algorithmic configuration depends on the instance at hand, are likely to produce better configurations. However, while PS approaches are generally problem-agnostic, PI ones require prior knowledge of the problem at hand, to efficiently encode each instance by a set of features.
PI approaches typically focus on learning a good surrogate map of the performance function, generally by performing regression: this approximation is used to direct the search in the configuration set. In \cite{ehmhh}, for example, linear basis function regression is used to approximate the target algorithm runtime, which is defined as a map of both features and configurations; then, for a new instance with known features, the learnt map is evaluated at all configuration points in an exhaustive search, to find the estimated best one. However, other approaches may be used: in \cite{lao}, for example, a map from instance features to optimal configuration is learnt by a neural network; in \cite{zarp_miqp} the ACP is restricted to a single binary parameter of CPLEX, and a classifier is then trained to predict it. In \cite{khal_branch}, instead, CPLEX is run on a given instance for a certain amount of computational resources, then a ranking ML model is trained, on-the-fly, to learn the ordering of branch and bound variables, and it is then used to predict the best branching variable, at each node, for the rest of the execution. Another approach, presented in \cite{isac}, first performs clustering on a set of instances, then uses the PS methodology described in \cite{gga} to find one good algorithmic configuration for each cluster. In \cite{VGC+19}, instead, instances are automatically clustered in the leaves of a trained decision tree, which also learns the best configuration for each leaf; at test time, a new instance is assigned to a leaf based on its features, and it receives the corresponding configuration.
The purpose of an ACP approach is to provide a good algorithmic configuration upon the arrival of an unseen instance. We call a methodology offline if the learning happens before that moment, which is the case for all the approaches cited above. Otherwise, we call an ACP methodology online; these approaches normally use reinforcement learning techniques (see, e.g.,~\cite{rlaoas}) or other heuristics \cite{reactive}.

\smallskip
In our approach, we define the performance of the target algorithm as a function of both features and controls, in order to account for the fact that the best configuration of a solver may vary among instances belonging to the same class of problems; this makes our approach PI. Moreover, we perform the PMLP only once, offline, which allows us to solve the resulting CSSP for any new instances in a matter of seconds.
What makes our approach stand out from other methodologies is that the learning phase is treated as white-box: the prediction problem of the PMLP is formulated as a MP, which conveniently allows the explicit embedding of a mathematical encoding of the estimated performance into the CSSP, as its objective/constraints. This is opposed to treating the learned predictor as a black-box, and therefore using it as an oracle in brute-force searches or similar heuristics, that typically do not scale as well as optimization techniques.

\section{The PMLP and the CSSP}\label{s:pmlp_cssp}

Let $\mathcal{A}$ be the target algorithm, and:
\begin{itemize}
	\item $\mathcal{C}_{\mathcal{A}}$ be the set of feasible configurations of $\mathcal{A}$. We assume that each configuration $c \in \mathcal{C}_{\mathcal{A}}$ can be encoded into a vector of binary and/or discrete values representing categorical and numerical parameters, and $\mathcal{C}_{\mathcal{A}}$ can be described by means of linear constraints;
	\item $\Pi$ be the problem to be solved, consisting of an infinite set of instances, and $\Pi' \subset \Pi$ be the (finite) set of instances used for the PMLP;
	\item $F_{\Pi}$ be the set of feature vectors used to describe instances, encoded by vectors of continuous or discrete/categorical values (in the latter case they are labelled by reals);
	\item $p_\mathcal{A} \,:\, F_{\Pi} \times \mathcal{C}_{\mathcal{A}} \longrightarrow \mathbb{R}$ be the performance function which maps a pair $(f,c)$ (instance feature vector, configuration) to the outcome of an execution of $\mathcal{A}$ (say in terms of the integrality gap reported by the solver after a time limit, but other measures are possible).
\end{itemize}
With the above definitions, the PMLP and the CSSP are detailed as follows.

\subsection{Performance Map Learning Phase}\label{ss:pmlp}

In the PMLP we use a supervised ML predictor, e.g.,~SVR, to learn the coefficient vector $\theta^\ast$ providing the parameters of a prediction model $\bar{p}_\mathcal{A}(\cdot,\cdot,\theta):F_{\Pi}\times\mathcal{C}_{\mathcal{A}}\rightarrow\mathbb{R}$ of the performance function $p_{\mathcal{A}}(\cdot,\cdot)$. The training set for the PMLP is
\begin{equation}
\mathcal{S} = \big\{(f_i, c_i, p_{\mathcal{A}}(f_i,c_i)) \;|\; i \in \{ 1 \dots s\} \big\} \ \subseteq\ F_{\Pi'} \times \mathcal{C}_{\mathcal{A}} \times \mathbb{R},
\label{eq:sample_set}
\end{equation}
where $s=|\mathcal{S}|$ and the training set labels $p_{\mathcal{A}}(f_i,c_i)$ are computed on the training vectors $(f_i,c_i)$. The vector $\theta^\ast$ is chosen as to hopefully provide a good estimate of $p_\mathcal{A}$ on points that do not belong to $\mathcal{S}$, with the details depending on the selected ML technology.

\subsection{Configuration Space Search Problem}\label{ss:cssp}

For a given instance $f$ and parameter vector $\theta$, $\mathsf{CSSP}(f,\theta)$ is the problem of finding the configuration with best estimated performance $\bar{p}_\mathcal{A}(f,c,\theta)$:
\begin{equation}
\mathsf{CSSP}(f,\theta) \equiv
\begin{array}{rrcl}
\min\limits_{c\in\mathcal{C}_{\mathcal{A}}} & \bar{p}_{\mathcal{A}}(f,c,\theta)
\end{array}. \label{eq:cssp}
\end{equation}
The actual implementation of $\mathsf{CSSP}(f, \theta)$ depends on the MP formulation selected to encode $\bar{p}_{\mathcal{A}}$, which may require auxiliary variables and constraints to define the properties of the ML predictor.
If $\bar{p}_{\mathcal{A}}$ yields an accurate estimate of $p_{\mathcal{A}}$, we expect the optimum $c^\ast_\mathsf{cssp}$ of $\mathsf{CSSP}(f,\theta)$ to be a good approximation of the true optimal configuration $c^\ast$ for solving $f$. However, we remark that a) $\mathsf{CSSP}(f,\theta)$ can be hard to solve, and b) it needs to be solved quickly (otherwise one might as well solve the instance $f$ directly).
Hence, incurring the additional computational overhead for solving the CSSP may be advantageous only when the instance at hand is ``hard''.
Achieving a balance between PMLP accuracy and CSSP cost is one of the challenges of this research.

\section{Experimental results}\label{s:experiments}

We tested our approach on 250 instances of the HUC problem and on 9 parameters of CPLEX, version 12.7. The PMLP and CSSP experiments were conducted on an Intel Xeon CPU E5-2620 v4 @ 2.10GHz architecture, while CPLEX was run on an Intel Xeon Gold 5118 CPU @ 2.30GHz. The pipeline was implemented in Python 3.6.8 \cite{VD09} and AMPL Version 20200110 \cite{ampl}. In the following, we detail the algorithmic set-up that we employed.

\subsection{Building the dataset}\label{ss:dataset}

\begin{enumerate}
	\item {\it Features}. The HUC is the problem of finding the optimal scheduling of a pump-storage hydro power station, where the commitment and the power generation of the plant must be decided in a short term period, in which inflows and electricity prices are previously forecast. The goal is to maximize the revenue given by power selling (see, e.g., \cite{shopa_huc}). The time horizon is fixed at 24 hours and the underlying hydro system is also fixed, so that all the instances have the same size. Thus, only 54 elements which vary from day to day are features: the date, 24 hourly prices, 24 hourly inflows, initial and target water volumes, upper and lower bound admitted on the water volumes. We encode them in a vector $f$ of 54 continuous/discrete components. All the instances have been randomly generated with an existing generator that accurately reproduces realistic settings.
	\item {\it Configuration parameters}. Thanks to preliminary tests, we select a subset of 9 discrete CPLEX parameters ({\tt fpheur}, {\tt dive}, {\tt probe}, {\tt heuristicfreq}, {\tt startalgorithm} and {\tt subalgorithm} from {\tt mip.strategy}; {\tt crossover} from {\tt barrier}; {\tt mircuts} and {\tt flowcovers}, from {\tt mip.cuts}), for each of which we consider between 2 and 4 different values. We then combine them so as to obtain 2304 parameter configurations. A configuration is encoded by a vector $c \in \{ \, 0 \,,\, 1 \, \}^{23}$, where each categorical parameter is represented by its incidence vector.
	\item {\it Performance measure}. We use the integrality gap to define $p_\mathcal{A}(f,c)$. It has been shown that MIP solvers can be affected by performance variability issues (see, e.g., \cite{perf_variability}), due to executing the solver on different computing platforms, permuting rows/columns of a model, adding valid but redundant constraints, performing apparently neutral changes to the solution process, etc. In order to tackle this issue, first we sample three different random seeds. For each instance feature vector $f$ and each configuration $c$, we then carry out the following procedure: (i) we run CPLEX (using the Python API) three times on the instance, using the different random seeds, for 60 seconds; (ii) we record the middle out of the three obtained performance values, to be assigned to the pair $(f,c)$.	At this point, our dataset contains $250\times2304=576000$ records. The performance measure thus obtained from CPLEX output, which we call $p_{\tt cpx}(f,c)$, usually contains some extremely large floating point values (e.g.,~whenever the CPLEX gap has a value close to zero in the denominator), which unduly bias the learning process. We deal with this issue as follows: we compute the maximum ${\bar{p}_{\tt cpx}}$, over all values of (the range of) $p_{\tt cpx}$, lower than a given threshold (set to $1\mathrm{e}{+5}$ in our experiments), re-set all values of $p_{\tt cpx}$ larger than the threshold to $\bar{p}_{\tt cpx}+100$, then rescale $p_{\tt cpx}$ so that it lies within the interval $[ \, 0 \,,\, 1 \, ]$. The resulting performance measure, which in the following we call $p_{\tt ml}(f,c)$, is also the chosen PMLP label.  Moreover, we solve each HUC instance and we record the value of its optimum; then, $\forall (f,c)$, we compute the {\it primal gap} $\varrho_\mathsf{prim}$ and the {\it dual gap} $\varrho_\mathsf{dual}$, i.e.: the distance between the optimal value of $f$ and the value of the feasible solution found, and the distance between the value of the tightest relaxation found and the optimal value (both over the optimum). We save $p_{\tt cpx}$, $p_{\tt ml}$, $\varrho_\mathsf{prim}$ and $\varrho_\mathsf{dual}$ in our dataset. We remark that setting the time-limit, imposed on CPLEX runs, to 60 seconds provides the solver enough time to move past the preliminary processing and to begin working on closing the gap, even for very hard instances (i.e., the ones with long pre-processing times); this allows us to measure the actual impact that different parameter configurations have on the chosen performance measure.
	\item {\it Feature engineering}. We process the date in order to extract the season, the week-day, the year-day, two flags called {\tt isHoliday} and {\tt isWeekend}, and we perform several sine/cosine encodings, that are customarily used to treat cyclical features. Moreover, we craft new features by computing statistics on the remaining 54 features. This task takes around 12 minutes to complete for the whole data set.
	\item {\it Splitting the dataset}. We randomly divide the instances into 187 In-Sample (IS) and 63 Out-of-Sample (OS), and split the dataset rows accordingly (430848 IS and 145152 OS). We use the IS data to perform Feature Selection (FS) and to train the SVR predictor; then, we assess the performance of the PMLP-CSSP pipeline both on OS instances, to test its generalization capabilities to unseen input, and on IS instances, to evaluate its performance on the data that we learn from, as detailed below.
	\item {\it Feature selection}. We use Python's \texttt{Pandas} \texttt{DataFrame}'s \texttt{corr} function to perform Pearson's Linear Correlation and \texttt{sklearn} \\ \texttt{RandomForestRegressor}'s \texttt{feature\_importances\_} attribute to perform decision trees' Feature Importance, in order to get insights on which features contribute the most to yield accurate predictions. A detailed explanation of the employed FS techniques falls outside of the scope of this document.
	In the following, we use the shorthand ``variables'' to refer to the whole list of columns of the learning dataset. In order to perform FS, we use a dedicated subset of the IS dataset, composed of 19388 records and only employed for this task; performing the selected FS techniques on this dataset takes around 8 minutes, and reduces $f$ to 22 components. For the configuration vectors we consider three FS scenarios: {\tt noFS}, {\tt kindFS} and {\tt aggFS}, yielding $c$ vectors with, respectively, 22, 14 and 10 components. We then filter the PMLP dataset according to the FS scenario at hand. However, after this filtering, the dataset may contain points with the same $(f,c)$ but different labels $p_{\tt ml}(f,c)$. Thus, for each instance: a) we delete the dataset columns that FS left out; b) we perform {\tt Pandas}'s {\tt group\_by} on the 22 columns chosen by FS, then c) compute the average $p_{\tt ml}$ of each group and use this as the new label (at this point, rows with the same $(f,c)$ have the same label); d) we remove the duplicate rows of the dataset by {\tt Pandas} {\tt drop\_duplicates}, keeping only one row. Lastly, we select $\sim$11200 points for the PMLP.
\end{enumerate}

\subsection{PMLP experimental setup} \label{ss:pmlp}

The PMLP methodology of choice in this paper is SVR. Its advantages are: (a) the PMLP for training an SVR can be formulated as a convex Quadratic Program (QP), which can be solved efficiently; (b) even complicated and possibly nonlinear performance functions can be learned by using the ``kernel trick'' \cite{svr}; (c) the solution of the PMLP for SVR provides a closed-form algebraic expression of the performance map $\bar{p}_{\mathcal{A}}$, which yields an easier formulation of $\mathsf{CSSP}(f,\theta)$. We use a Gaussian kernel during SVR training, which is the default choice in absence of any other meaningful prior \cite{svm_and_kernels}. We assess the prediction error of the predictor by Nested Cross Validation (NCV) \cite{ncv}; furthermore, our training includes a phase for determining and saving the hyperparameters and the model coefficients of the SVR. These two tasks take approximately 4 hours in the {\tt aggFS} and in the {\tt kindFS} scenarios, and 5 hours in the {\tt noFS} one. We use Python's \texttt{sklearn.model\_selection.RandomizedSearchCV} for the inner loop of the NCV and a customized implementation for the outer loop, and  \texttt{sklearn.svm.SVR} as the implementation of choice for the ML model.

A common issue in data-driven optimization is that using customary ML error metrics may not lead to good solutions of the optimization problem (see, for example, \cite{pr+opt}). We tackled this issue by comparing the following metrics for the CV-based hyperparameter tuning phase, both computed on $p_{\tt ml}$: the classical Mean Absolute Error $\mathsf{MAE} = \sum_{i \in S} |p_i - \bar{p}_i|$, where $p_i = p_{\tt ml}(f_i,c_i)$ and $\bar{p}_i = \bar{p}_\mathcal{A}(f_i,c_i)$; the custom metric $ \mathsf{cMAE}_\delta = \sum_{i \leq s} L_\delta(p_i, \bar{p}_i)$, $\delta \in \{0.2, 0.3, 0.4\}$, where
\begin{equation*}
L_\delta(p_i, \bar{p}_i) =
\left\{\begin{array}{ll}
(\bar{p}_i - p_i) \cdot (1 + \frac{1}{1 + exp(p_i - \bar{p}_i)}) \quad & \text{if } p_i \leq \delta \mbox{ and } \bar{p}_i > p_i \\
(p_i - \bar{p}_i) \cdot (1 + \frac{1}{1 + exp(\bar{p}_i - p_i)}) \quad & \text{if } p_i \geq 1 - \delta \mbox{ and } \bar{p}_i < p_i \\
(p_i - \bar{p}_i) \quad & \text{if } \delta \leq p_i \leq 1 - \delta\\
0 \quad & \text{otherwise.}
\end{array}\right.
\end{equation*}

\subsection{CSSP experimental setup}\label{ss:cssp}

The choice of a Gaussian kernel in the SVR formulation makes the CSSP a MINLP with a nonconvex objective function $\bar{p}_\mathcal{A}$. More precisely, for a given instance with features $\bar{f}$, our CSSP is:
\begin{equation}
\textstyle
\min\limits_{c\in\mathcal{C}_{\mathcal{A}}}\; \sum_{i=1}^s \alpha_i\,
\exp\Big({\,-\gamma \|(f_i, c_i) - (\bar{f}, c)\|_2^2} \Big)
\label{eq:cssp_svr_obj}
\end{equation}
where, for all $i\le s$, $(f_i, c_i)$ belong to the training set, $\alpha_i$ are the dual solutions of the SVR, $\gamma$ is the scaling parameter of the Gaussian kernel, and $\mathcal{C}_{\mathcal{A}}$ is defined by mixed-integer linear programming constraints encoding the dependences/compatibility of the configurations. We use AMPL to formulate the CSSP, and the nonlinear solver Bonmin \cite{bonmin_algo}, manually configured (with settings \texttt{heuristic\_dive\_fractional yes}, \texttt{algorithm B-Hyb}, \texttt{heuristic\_feasibility\\\_pump yes}) and with a time limit of 60 seconds, to solve it; then we retrieve, for each instance $f$, the Bonmin solution $c^\ast_{\mathsf{bm}}$. Since we have enumerated all possible configurations, we can also compute the ``true'' global optimum $c^\ast_{\mathsf{cssp}}$ $= \arg\min\{\bar{p}_\mathcal{A}(f,c), c \in \mathcal{C}_\mathcal{A}\}$ for sake of comparison.  In Table~\ref{t:cssp_quality_barp_FS}, we report the percentage of cases where $c^\ast_{\mathsf{bm}}$ $=$ $c^\ast_{\mathsf{cssp}}$ (``\%glob.~mins'') and, for all the instances where this is not true, the average distance between $\bar{p}_\mathcal{A}(c^\ast_{\mathsf{bm}})$ and $\bar{p}_\mathcal{A}(c^\ast_{\mathsf{cssp}})$, over all the instances of the considered set (``avg loc.~mins''); we also report the average CSSP solution time (``CSSP time'').
\begin{table}[ht]
\centering
\scalebox{0.8}{
	\begin{tabular}{c|c|c|c|c}
		set   & FS & \%glob. mins & avg loc. mins & CSSP time \\
		\hline
		\multirow{3}[1]{*}{IS} & {\tt noFS}  & 83.69 & $\bm{2.013\mathrm{e}{-02}}$ & 15.36 \\
		& {\tt kindFS} & {\bf 87.70 }& 2.758$\mathrm{e}{-02}$ & 10.73 \\
		& {\tt aggFS} & 84.49 & 3.122$\mathrm{e}{-02}$ & 5.95 \\
		\hline
		\multirow{3}[0]{*}{OS} & {\tt noFS}  & 85.32 & $\bm{1.594\mathrm{e}{-02}}$ & 13.44 \\
		& {\tt kindFS} & 90.48 & 2.014$\mathrm{e}{-02}$ & 12.23 \\
		& {\tt aggFS} & {\bf 93.25} & 1.957$\mathrm{e}{-02}$ & 6.15 \\
	\end{tabular}%
}
\vspace{2mm}
\caption{Quality of Bonmin's solutions, w.r.t.~$\bar{p}$}
\label{t:cssp_quality_barp_FS}
\end{table}
The {\tt kindFS} and the {\tt aggFS} scenarios achieve better results, in terms of ``\%glob mins'', than the {\tt noFS} one.
Furthermore, the local optima found by Bonmin are quite good ones: they are never larger than 3.2\% (see ``avg loc.~mins'') of $c^\ast_\mathsf{cssp}$.
The time that Bonmin takes to solve the CSSP is reduced from the {\tt noFS} scenario to the {\tt aggFS} one. This is due to the fact that the {\tt kindFS} and the {\tt aggFS} CSSP formulations have less variables than the {\tt noFS} one, and so they are easier to solve. Bonmin needs, on average, less than 16 seconds to solve any CSSP; however, devising more efficient techniques to solve the CSSP (say, reformulations, decomposition, \dots) might be necessary if our approach is scaled to considerably more algorithmic parameters.

\subsection{Results}\label{ss:pipeline}
In order to assess the performance of the approach, we retrieve $p_{\tt cpx}(c^\ast_{\mathsf{bm}})$, $p_{\tt cpx}(c_{\mathsf{cpx}})$ (CPLEX default configuration), the primal and dual gap of $c^\ast_{\mathsf{bm}}$ and $c_{\mathsf{cpx}}$ from the filtered dataset, for every IS and OS instance.
\begin{table}[h]
	\centering
	\scalebox{0.8}{
		\begin{tabular}{c|c|ccc|ccc}
			set   & FS    & \%w   & \%w+d & \%w$_{\mbox{nond}}$ & avg d & avg w & avg l \\
			\hline
			\multirow{6}[-10]{*}{IS} & {\tt noFS } & 47.06 & 96.12 & 92.39 & 0 & 3.400$\mathrm{e}{+14}$ & $\bm{6.493\mathrm{e}{+19}}$ \\
			& {\tt kindFS} & 48.91 & 98.02 & {\bf 96.12} & 0 & 3.415$\mathrm{e}{+14}$ & 8.895$\mathrm{e}{+19}$ \\
			& {\tt aggFS} & {\bf 48.98} & 97.79 & 95.70 & 0 & $\bm{3.438\mathrm{e}{+14}}$ & 8.530$\mathrm{e}{+19}$ \\
			\hline
			\multirow{6}[-10]{*}{OS} & {\tt noFS } & 33.73 & 83.73 & 67.50 & 0 & 3.148$\mathrm{e}{+14}$ & 7.042$\mathrm{e}{+19}$ \\
			& {\tt kindFS} & {\bf 36.64} & 82.61 & {\bf 68.06} & 0 & 3.051$\mathrm{e}{+14}$ & 5.385$\mathrm{e}{+19}$ \\
			& {\tt aggFS} & 33.86 & 77.05 & 59.78 & 0 & $\bm{3.335\mathrm{e}{+14}}$ & $\bm{4.486\mathrm{e}{+19}}$ \\
		\end{tabular}%
	}
	\vspace{2mm}
	\caption{Pipeline quality w.r.t.~$p_{\tt cpx}(c^\ast_{\mathsf{bm}})$, by FS scenario and IS/OS set}
	\label{t:pipeline_quality_bonmin_FS}
\end{table}

Tab.~\ref{t:pipeline_quality_bonmin_FS} shows:
the wins ``\%w'' and the non-worsenings ``\%w+d'', i.e., the percentage of instances such that $p_\mathsf{cpx}(c^\ast_\mathsf{bm})$ is $<$ or $\leq$ than $p_\mathsf{cpx}(c_\mathsf{cpx})$, by the first sixteen decimal digits of $p_\mathsf{cpx}$, in scientific notation; the wins-over-nondraws ``\%w$_{\mbox{nond}}$'', i.e., the percent wins over the instances such that $p_\mathsf{cpx}(c^\ast_\mathsf{bm}) \neq p_\mathsf{cpx}(c_\mathsf{cpx})$; the average $|p_\mathsf{cpx}(c_\mathsf{cpx}) - p_\mathsf{cpx}(c^\ast_\mathsf{bm})|$, over all the instances which score a win (``avg w'') or a loss (``avg l''); the average $|p_\mathsf{cpx}(c^\ast_\mathsf{bm}) - \underline{p}_\mathsf{cpx}|$ over all the other instances, where $\underline{p}_\mathsf{cpx} = \min_{c \in \mathcal{C}_\mathcal{A}} p_\mathsf{cpx}(f, c)$ for a given $f$ (``avg d'').
The ``\%w'', ``\%w+d'' and ``\%w$_{\mbox{nond}}$'' are higher on IS instances than on the OS ones. The IS instances are used as the training set, so it is not surprising that $\bar{p}_\mathcal{A}$ is less accurate at OS instances; this results in worse CSSP solutions for OS instances.
The fact that ``avg d'' is always 0 implies that, whenever CPLEX --- configured by $c_\mathsf{cpx}$ --- can close the gap, our $c^\ast_\mathsf{bm}$ proves to be just as efficient. In the nondraws, our approach shows consistent gains with respect to $c_\mathsf{cpx}$, both on IS and OS instances. From this we gather that $\bar{p}_\mathcal{A}$ provides an accurate approximation of $p_\mathsf{cpx}$'s global minima, even at points outside the training set.
The {\tt noFS} scenario presents the worst results. The {\tt aggFS} scenario achieves the highest percent wins on IS instances; it also provides the largest average wins and the smallest average loss, overall. However, the {\tt kindFS} scenario presents the best performance, by the highest ``\%w'' on OS instances and by the best overall ``\%w$_{\mbox{nond}}$''.
\begin{table}[ht]
	\centering
	\scalebox{0.7}{
		\begin{tabular}{c|c|c|c|c|c|c|c|c}
			set   & FS    & metric & \%feas $c^\ast_\mathsf{bm}$ & \%feas $c_\mathsf{cpx}$ & $\varrho_\mathsf{prim}^\mathsf{bm}$ & $\varrho_\mathsf{prim}^\mathsf{cpx}$ & $\varrho_\mathsf{dual}^\mathsf{bm}$ & $\varrho_\mathsf{dual}^\mathsf{cpx}$  \\
			\hline
			\multirow{12}[6]{*}{IS} & \multirow{4}[2]{*}{noFS} & $\mathsf{cMAE}_{.2}$ & 96.79 & \multirow{12}[6]{*}{100} & 1.504$\mathrm{e}{-01}$ & \multirow{12}[6]{*}{5.812$\mathrm{e}{-01}$} & 9.708$\mathrm{e}{-02}$ & \multirow{12}[6]{*}{1.018$\mathrm{e}{-01}$}  \\
			&       & $\mathsf{cMAE}_{.3}$ & 95.72 &       & 1.193$\mathrm{e}{-01}$ &       & 8.104$\mathrm{e}{-02}$ &  \\
			&       & $\mathsf{cMAE}_{.4}$ & 98.40 &       & 1.143$\mathrm{e}{-01}$ &       & 7.253$\mathrm{e}{-02}$ &  \\
			&       & $\mathsf{MAE}$   & \textbf{98.93} &       & $\bm{1.058\mathrm{e}{-01}}$ &       & $\bm{4.698\mathrm{e}{-02}}$ &   \\
			\cline{2-4}\cline{6-6}\cline{8-8}      & \multirow{4}[2]{*}{{\tt kindFS}} & $\mathsf{cMAE}_{.2}$ & 98.75 &       & $\bm{6.975\mathrm{e}{-02}}$ &       & $\bm{3.792\mathrm{e}{-02}}$ &   \\
			&       & $\mathsf{cMAE}_{.3}$ & \textbf{99.38} &       & 8.271$\mathrm{e}{-02}$ &       & 4.245$\mathrm{e}{-02}$ &  \\
			&       & $\mathsf{cMAE}_{.4}$ & 98.48 &       & 7.085$\mathrm{e}{-02}$ &       & 5.015$\mathrm{e}{-02}$ &  \\
			&       & $\mathsf{MAE}$   & 96.43 &       & 7.183$\mathrm{e}{-02}$ &       & 4.159$\mathrm{e}{-02}$ &   \\
			\cline{2-4}\cline{6-6}\cline{8-8}      & \multirow{4}[2]{*}{{\tt aggFS}} & $\mathsf{cMAE}_{.2}$ & \textbf{98.04} &       & 7.850$\mathrm{e}{-02}$ &       & 5.179$\mathrm{e}{-02}$ &   \\
			&       & $\mathsf{cMAE}_{.3}$ & 97.86 &       & $\bm{7.669\mathrm{e}{-02}}$ &       & 5.111$\mathrm{e}{-02}$ &  \\
			&       & $\mathsf{cMAE}_{.4}$ & 97.59 &       & 9.267$\mathrm{e}{-02}$ &       & 6.680$\mathrm{e}{-02}$ &  \\
			&       & $\mathsf{MAE}$   & 97.59 &       & 7.747$\mathrm{e}{-02}$ &       & $\bm{4.998\mathrm{e}{-02}}$ &   \\
			\hline
			\multirow{12}[5]{*}{OS} & \multirow{4}[2]{*}{{\tt noFS}} & $\mathsf{cMAE}_{.2}$ & 87.30 & \multirow{12}[5]{*}{100} & $\bm{6.968\mathrm{e}{-02}}$ & \multirow{12}[5]{*}{4.633$\mathrm{e}{-01}$} & 7.931$\mathrm{e}{-02}$ & \multirow{12}[5]{*}{8.121$\mathrm{e}{-02}$}  \\
			&       & $\mathsf{cMAE}_{.3}$ & 92.06 &       & 1.231$\mathrm{e}{-01}$ &       & 8.500$\mathrm{e}{-02}$ &  \\
			&       & $\mathsf{cMAE}_{.4}$ & 88.89 &       & 1.305$\mathrm{e}{-01}$ &       & 9.687$\mathrm{e}{-02}$ &  \\
			&       & $\mathsf{MAE}$   & \textbf{93.65} &       & 1.371$\mathrm{e}{-01}$ &       & $\bm{6.611\mathrm{e}{-02}}$ &   \\
			\cline{2-4}\cline{6-6}\cline{8-8}      & \multirow{4}[2]{*}{{\tt kindFS}} & $\mathsf{cMAE}_{.2}$ & 89.68 &       & 1.471$\mathrm{e}{-01}$ &       & $\bm{6.471\mathrm{e}{-02}}$ &   \\
			&       & $\mathsf{cMAE}_{.3}$ & 87.83 &       & $\bm{1.261\mathrm{e}{-01}}$ &       & 1.189$\mathrm{e}{-01}$ &  \\
			&       & $\mathsf{cMAE}_{.4}$ & 91.01 &       & 1.541$\mathrm{e}{-01}$ &       & 1.283$\mathrm{e}{-01}$ &  \\
			&       & $\mathsf{MAE}$   & \textbf{92.86} &       & 1.486$\mathrm{e}{-01}$ &       & 1.078$\mathrm{e}{-01}$ &   \\
			\cline{2-4}\cline{6-6}\cline{8-8}      & \multirow{4}[1]{*}{{\tt aggFS}} & $\mathsf{cMAE}_{.2}$ & \textbf{89.42} &       & 2.054$\mathrm{e}{-01}$ &       & 9.062$\mathrm{e}{-02}$ &   \\
			&       & $\mathsf{cMAE}_{.3}$ & 88.89 &       & 2.061$\mathrm{e}{-01}$ &       & 9.395$\mathrm{e}{-02}$ &  \\
			&       & $\mathsf{cMAE}_{.4}$ & 88.10 &       & 2.075$\mathrm{e}{-01}$ &       & 9.874$\mathrm{e}{-02}$ &  \\
			&       & $\mathsf{MAE}$   & 88.62 &       & $\bm{1.792\mathrm{e}{-01}}$ &       & $\bm{7.733\mathrm{e}{-02}}$ &  \\
		\end{tabular}%
	}
	\vspace{2mm}
	\caption{Quality of the solutions attained by CPLEX, configured by $c^\ast_{\mathsf{bm}}$ and $c_{\mathsf{cpx}}$, solving HUC instances, aggregated by IS/OS set, FS scenario and PMLP metric}
	\label{t:pipeline_quality_bonmin_feas}
\end{table}
In Tab.~\ref{t:pipeline_quality_bonmin_feas}, we report the percentage of IS and OS instances such that CPLEX, configured by $c^\ast_\mathsf{bm}$ and $c_\mathsf{cpx}$, manages to find a feasible solution (``\%feas $c^\ast_\mathsf{bm}$'' and ``\%feas $c_\mathsf{cpx}$''); for those instances, the columns ``$\varrho_\mathsf{prim}^\mathsf{bm}$'',  ``$\varrho_\mathsf{prim}^\mathsf{cpx}$'', ``$\varrho_\mathsf{dual}^\mathsf{bm}$'' and ``$\varrho_\mathsf{dual}^\mathsf{cpx}$'' report the average primal and dual gap achieved by the solver. CPLEX's default configuration always allows the solver to obtain a feasible solution within the time-limit. Our approach presents similar results on IS instances (``\%feas $c^\ast_\mathsf{bm}$'' is approximately 98\%), but it has slightly worse performances on the OS ones (``\%feas $c^\ast_\mathsf{bm}$'' is around 90\%). However, the primal/dual gaps provided by our methodology are always better than those achieved by using CPLEX's default setting; actually, they are up to an order of magnitude smaller,  in IS instances. The $\mathsf{cMAE_{.4}}$ and the $\mathsf{MAE}$ metrics provide the highest ``\%feas $c^\ast_\mathsf{bm}$'', respectively, on IS instances (above 98\%) and OS instances (around 92\%). On average, the $\mathsf{MAE}$ also provides the best primal/dual gaps for IS instances, while the $\mathsf{cMAE_{.4}}$ is the best choice for OS instances. Lastly, while the {\tt kindFS} scenario prevails on IS instances, the {\tt noFS} one dominates on the OS ones.

\section{Conclusions}
The methodology presented in this paper conflates ML and MP techniques to solve the ACP. All in all, the results show that the approach is promising, in that the configurations that it provides typically have better primal/dual gaps than CPLEX’s. The fact that, in a small fraction of the cases, no feasible solution is found depends on how infeasibility is encoded by $p_{\tt ml}$; an easy fix to this issue would be to use a performance measure promoting feasibility. Notably, it is interesting that, in some cases, the custom $\mathsf{cMAE}$ error metric outperforms the classical $\mathsf{MAE}$. We also observed that using FS techniques is conducive to much easier CSSP formulations, without overly affecting the quality of the solutions. Overall, since choices taken at any point in the pipeline affect its final outcome, a number of details have to be carefully considered for the approach to work.

\ifarxiv
\bibliographystyle{unsrt}
\else
\bibliographystyle{splncs04}
\fi
{\bibliography{biblio_giommazz}}

\clearpage
\appendix
\onecolumn
\end{document}